\documentclass[11pt]{ensmath}

\EnsMath pages 1-21

\usepackage{url}
\usepackage{amssymb, latexsym, amsfonts, amscd, amsbsy}

\newtheorem{Thm}{Theorem}
\newtheorem{Conj}{Conjecture}

\newtheorem{Prop}{Proposition}

\title[Two-dimensional lattices with few distances]{Two-dimensional lattices with few distances} 
 
\author[Pieter Moree \and  Robert Osburn]{Pieter \sn{Moree} \and Robert \sn{Osburn}}
 
\begin{document}
 
\maketitle

\begin{abstract}
We prove that of all two-dimensional lattices of 
covolume 1 the hexagonal lattice 
has asymptotically
the fewest distances. An analogous result for dimensions 3 to 8 was
proved in 1991 by Conway and Sloane. Moreover, we give a survey of some related
literature, in particular progress on a conjecture from 1995 due to Schmutz
Schaller. 
\end{abstract}

\section{Introduction}
It is an old problem in combinatorial geometry how to place a given number
of distinct points in $n$-dimensional Euclidean space so as to minimize
the total number of distances they determine. Conway and Sloane \cite{CS}
conjecture that, for all $N$ sufficiently
large, the optimal set of $N$ points in $n$-dimensional space
will be a subset of an $n$-dimensional lattice having minimal {\it Erd{\H o}s number}. 
In real Euclidean space $\mathbb R^n$ equipped with inner product $(v,w)=v\cdot w$, 
a lattice $L$ consists of all integral linear combinations
$$v=\lambda_1 v_1+\cdots+\lambda_n v_n,~\lambda_i\in \mathbb Z,$$
of $n$ linearly independent vectors $v_1,\ldots,v_n$. The vectors $v_1,\ldots,v_n$ form
an integral basis for $L$, and 
$$f(\lambda)=(v,v)=\lambda A \lambda^{\rm tr},~\lambda=(\lambda_1,\ldots,\lambda_n),~A=(a_{ij}),~
a_{ij}=(v_i,v_j),$$
is the corresponding quadratic form. The various integral bases for $L$ yield integrally
equivalent quadratic forms.
Suppose $n\ge 2$.
The
Erd{\H o}s number of an $n$-dimensional lattice $L$ is given by
\begin{equation}
\label{allebegin}
E_L=F_Ld^{1/n},
\end{equation}
where $d$ is the determinant of the lattice and $F_L$, its {\it population
fraction}, is given by
$$F_L=\lim_{x\rightarrow \infty}{N_L(x)\sqrt{\log x}\over x}{\rm ~if~}n=2,~
F_L=\lim_{x\rightarrow \infty}{N_L(x)\over x}{\rm ~if~}n\ge 3,$$
where $N_L(x)$ is the {\it population function} associated to the corresponding
quadratic form, i.e., the number of values not exceeding $x$ taken by
the form. The Erd{\H o}s number is the population fraction when the lattice
is normalized to have covolume 1. Conway and Sloane \cite{CS} proved that 
for $n\ge 3$ the lattices with minimal Erd{\H o}s number are (up to a scale factor) the
even lattices of minimal determinant. For $2\le n\le 10$ the even lattices 
of minimal determinant are
unique:

\begin{equation}
\label{ADE}
A_2,~A_3\cong D_3,~D_4,~D_5,~E_6,~E_7,~E_8,~E_8\oplus A_1,~E_8\oplus A_2.
\end{equation}

Actually Conway and Sloane also claimed the result for $n=2$, relying on a preprint (in 1990) of Warren D. Smith \cite{Smith}. However, the preprint was never published and
this induced Schmutz Schaller \cite[p. 200]{Schaller1} to write `the case $n=2$
seems to be open'. It is the purpose of this paper to dispose of this case (in Theorem 
\ref{main})
and thus to `complete' the 
Conway and Sloane result. In doing so, we have made use of results that have become
available only very recently. In particular, we use an explicit formula for the number of genera of discriminant $D$ representing a positive integer $n$ (see Theorem \ref{sunny}) and an improved lower bound on the Euler phi function $\varphi(n)$ for $n$ odd (see (\ref{nicolas})). \\ 
\indent Let $\Sigma$ denote the hexagonal lattice of
covolume 1, that is,
$$\Sigma= {1\over \sqrt{3}} A_2 =  \sqrt{2\over \sqrt{3}}\left(\mathbb Z
\left({1\atop 0}\right)\oplus \mathbb Z\left({1/2\atop 
\sqrt{3}/2}\right)\right).$$
The associated quadratic form is $(X^2+XY+Y^2)2/\sqrt{3}$.

\begin{Thm}
\label{main}
If $L$ is any two-dimensional lattice not isometric to $\Sigma$, then $E_L$, the
Erd{\H o}s number of $L$, satisfies 
\begin{equation}
\label{basaleongelijkheid}
E_L>E_{\Sigma}=2^{-3/2}3^{1/4}\prod_{p\equiv 2({\rm mod~}3)}
{1\over \sqrt{1-{1/p^2}}}=0.553311775832479\cdots
\end{equation}
In other words, of all the two dimensional lattices of covolume 1, $\Sigma$ has
asymptotically the fewest distances. Moreover, given any real number $r$ the set
of non-homothetic lattices $L$ 
such that $E_L<r$ is finite and can be explicitly determined.
\end{Thm}

In fact it turns out furthermore that if $E_L$ is finite, then there is 
a homothetic lattice $L'$ such that $E_L=E_{L'}$ and the
quadratic form associated to $L'$ has integer coefficients and is
primitive (for $\Sigma$ this is $X^2+XY+Y^2$). Moreover, 
$E_L$ only will depend on the discriminant $D$ of the associated 
quadratic form. To stress this, we write $E(D)$ rather than $E_L$.

\subsection{On a conjecture of Schmutz Schaller}
In \cite[p. 20]{Schaller1} Schmutz Schaller, motivated by considerations from
hyperbolic geometry, proposed for dimensions 2 to 8 a daring strengthening of 
Theorem \ref{main} and (part of) the Conway and Sloane result: 

\begin{Conj}
\label{conj1}
In dimensions $2$ to $8$ the even lattices
with minimal determinant have `maximal lengths', meaning that their
length spectrum dominates the length spectrum of every other lattice
of the same dimension and covolume at every position.
\end{Conj}

Schmutz Schaller \cite{Schaller0} proved an analogue of this conjecture in the hyperbolic
case.
Given any lattice $L$ one can define the sequence $0<d_1<d_2<...$ of distances
between lattice points that occur in this lattice (the {\it length spectrum}). (It
is very important that in this definition we do not care about the multiplicities 
of these lengths.) The number $d_k$ is called the {\it $k$-th length} of $L$.
Given any other length spectrum $0<l_1<l_2<\ldots$ we say that the
former length spectrum {\it totally dominates} the latter if $d_i\ge l_i$ for every $i\ge 1$.
This can be reformulated in terms of $N_L(x)$: the length spectrum $L_1$ {\it totally dominates} that of $L_2$ iff $N_{L_1}(x)\le N_{L_2}(x)$
for every $x>0$.\\
\indent Let $S=\mathbb Z[i]$ denote the square lattice and
$H=\mathbb Z[\zeta_3]$ the hexagonal lattice. Schmutz Schaller \cite{Schaller1}
conjectured
that the hexagonal length spectrum should dominate that of the square lattice, 
that is he conjectured that $N_{H}(x)\le N_{S}(x)$ for every $x>0$, to make the
point that even a partial version of his conjecture should be difficult
to establish. Indeed, the first author and te Riele \cite{MR}, refining 
techniques from \cite{Mchebbie}, managed to prove this only 
after considerable effort, also numerical effort. Their approach, however, does not seem to
offer any hope of establishing the general conjecture.\\
\indent {}From the work of Korkine and Zolotareff (in the 19th century) and Blichfeldt 
(cf. \cite[Chapter 9]{SO} and \cite{Blik}) it follows that Conjecture \ref{conj1} is
true in the $1$-length case, i.e., the lattices in (\ref{ADE}) have maximal minimal
positive length amongst those of the same dimension, after scaling to the same covolume. For 
a list of these lengths see, e.g. \cite[p. 204]{SO}.\\
\indent A two-dimensional lattice is said to be {\it arithmetic} iff there exists a real number $\lambda$
such that $\lambda L$ is isometric to a $\mathbb Z$-submodule of rank two in an imaginary
quadratic number field, otherwise it is said to be {\it non-arithmetic}. 
K\"uhnlein \cite{Kuhnlein} proved that a two-dimensional lattice is arithmetic iff there
are at least 3 pairwise linearly independent vectors in it having the same length.
As a consequence it is easy to show that
$N_L(x)\sim c(L)x$ for some positive constant $c(L)$ in case 
$L$ is non-arithmetic. It follows from this
that a non-arithmetic lattice does not have a finite Erd{\H o}s number. 
K\"uhnlein \cite{Kuhnlein} proved furthermore that the length spectrum of $\Sigma$ 
totally dominates the length spectrum
of every non-arithmetic lattice of covolume 1. Thus in
order to prove Conjecture \ref{conj1} for dimension 2 it suffices to prove
that the length spectrum of $\Sigma$ totally dominates the length spectrum of
every arithmetic lattice of covolume 1.

\section{Population fraction of binary quadratic forms}
Let $f(X,Y)=aX^2 +bXY + cY^2$ be a  positive definite binary quadratic form 
with discriminant 
$D_f=b^2-4ac$ and $a,b$ and $c$ real numbers.
Let $B_f(x)$ count the number of positive real numbers $r\le x$ that can be represented by $f$.

In the course of history the problem of estimating $B_f(x)$ has attracted considerable interest.
A classical result of 
Landau \cite{L} states that, as $x$ tends to infinity,
$$B_{f_1}(x) \sim C(f_1){x\over \sqrt{\log x }},$$
where $C(f_1)$ is an explicit constant and $f_1(X,Y)=X^2+Y^2$.
Precisely, $C(f_1)$ is of the form
$$
C(f_1)=\frac{1}{\sqrt{2}} \prod_{p\equiv 3 ({\rm mod~}4)} (1-p^{-2})^{-1/2}.
$$
Note that $B_{f_1}(x)=N_{S}(x)$. 

A similar result was claimed by Srinivasa Ramanujan in his celebrated first
letter to Hardy (written in 1912), cf. \cite{MoCa}. The constant $C(f_1)$ is 
now called the {\it Landau-Ramanujan constant}, cf. \cite[Section 2.3]{Finch}.
Ramanujan even claimed that it ought to be true that

\begin{equation}
\label{een} 
N_S(x)=C(f_1)\int_2^x{dt\over \sqrt{\log t}}+O(x^{1/2+\epsilon}).
\end{equation}

Note the analogy with the prime number theorem under assumption of the
Riemann Hypothesis. 
This states that $\pi(x)$, the number of primes $p\le x$, satisfies $\pi(x)=\int_2^x{dt/\log t}
+O(x^{1/2+\epsilon})$, on assumption of the Riemann Hypothesis. 
It was folklore that Landau's method 
could be easily adapted to show that $N_S(x)$ satisfies an asymptotic series expansion in the 
sense of Poincar\'e:

\begin{equation}
N_S(x)=C(f_1){x\over \sqrt{\log x}}\left(1+{r_1\over \log x}+{r_2\over \log ^2 x}+
\cdots+{r_n\over \log ^m x}+O\left({1\over \log^{m+1}x}\right)\right),
\end{equation}

\noindent where $m\ge 1$ is an arbitrary integer.
A proof of this was finally written down by J.-P. Serre \cite{Serre} for the larger
class of so called Frobenian multiplicative functions. Note that Ramanujan's conjecture
implies, by partial integration of the main term, that
$$N_S(x)=C(f_1){x\over \sqrt{\log x}}\left(1+{s_1\over \log x}+{s_2\over \log ^2 x}+
\cdots+{s_m\over \log ^m x}+O\left({1\over \log^{m+1}x}\right)\right),$$
with $s_j=(2j-1)!/((j-1)!2^{2j-1})$ and $m\ge 1$ an arbitrary integer.
Ramanujan's conjecture was shown to be false by Shanks \cite{Shanks} who 
proved that $s_1\ne r_1$. In a celebrated unpublished (during
his lifetime)
paper on the partition and tau function \cite{bono}, Ramanujan made conjectures similar to
(\ref{een}) concerning the divisiblity of the Ramanujan tau function by certain special
primes. These conjectures were all shown to be false by the first author \cite{moreetau}. However, 
Rankin had shown earlier that asymptotically these conjectures are correct.\\
\indent Paul Bernays (of later fame in logic and for many
years assistant to Hilbert \cite{Reid}) was a PhD student of Landau's at G{\"o}ttingen. In his 
1912 thesis Bernays \cite{B} studied the question 
of finding an asymptotic formula similar to that of Landau's, but now in case $f$ 
is a primitive positive definite binary quadratic form having negative discriminant $D_f$.
Bernays' proved that, as $x$ tends to infinity,

\begin{equation}
\label{R.8}
B_f(x)=C(f){x\over \sqrt{\log x}}+O\left({x\over (\log x)^{1/2+\delta}}\right),
\end{equation}

\noindent where the constant $C(f)$ is positive and depends only on the discriminant $D_f$ of 
$f$ and $\delta<{\rm min}(1/h,1/4)$, where $h$ denotes the number of reduced
quadratic forms having the same discriminant as $f$. It turns out that the dependence of 
$C(f)$ on $D_f$ is not very strong; 
$C(f)=D_f^{o(1)}$.\\
\indent Bernays' result allows various generalisations: one could ask for simultaneous
representation of $n$ by various quadratic forms or by norm forms. A lot of work
in this direction was carried out by Odoni, cf. \cite{Odoni-2, Odoni-1}. Blomer recently pointed 
out that Bernays' method can be used to disprove 
a conjecture of Erd{\H o}s. The falsity of this conjecture was claimed
earlier by Odoni \cite{Odoni}, but his paper seems to contain some
obscurities. Erd{\H o}s conjectured that the number $V(x)$ of integers
not exceeding $x$ that are sums of two squareful integers satisfies
$V(x)\asymp x/\sqrt{\log x}$, where an integer $n$ is called squareful
if $p|n$ implies that $p^2|n$ for all primes $p$. Since every squareful
integer $n$ can uniquely be written as $n=a^3b^2$ with $\mu(a)\ne 0$, one
can write
$$V(x)=\#\{1\le n\le x:\exists {\bf a}=(a_1,a_2)\in \mathbb N^2:a_1^3X^2+a_2^3Y^2
{\rm ~represents~}n\}.$$
Thus one can estimate $V(x)$ if one can deal with $B_f(x)$ with some
uniformity in $f$ (or rather the discriminant of $f$). In Bernays' method
the dependence on $D$ can be made explicit. This yields
$B_f(x)\gg _{\epsilon} |D|^{-\epsilon}x/\sqrt{\log x}$ uniformly at least
in $D=O((\log \log x)^{1/2})$. This result can be used to show that
Erd{\H o}s' conjecture is false. By a more refined method Blomer 
\cite{Blomer, Blomer1}
even showed that $V(x)=x(\log x)^{-\alpha+\epsilon}$, where
$\alpha=1-2^{-1/3}=0.206\cdots$. Moreover, Blomer and Granville \cite{BG} conjecture 
that $V(x)\asymp x(\log \log x)^{2^{2/3}-1}(\log x)^{2^{-1/3}-1}$ and prove the 
upper bound, failing to obtain the conjectured lower bound only by a power of
$\log \log x$.\\
\indent Bernays' result can be used to infer the following alternative characterisation of
arithmetic lattices.

\begin{Prop}
\label{prop1}
A two-dimensional lattice has a finite Erd{\H o}s number iff it is arithmetic.
\end{Prop}

\begin{proof} 
We have already seen that a non-arithmetic lattice does not have a finite 
Erd{\H o}s number. If the 
lattice is arithmetic then, possibly after scaling, the associated quadratic
form has integer coefficients. The result then follows from Bernays' theorem and
the definition (\ref{allebegin}) for $n=2$. \qed
\end{proof}

We say that the quadratic form $f=[a,b,c]$ is projectively equivalent with
$g=[a',b',c']$ if the vectors $(a,b,c)$ and $(a',b',c')$ are projectively equivalent.
If $g$ is projectively equivalent to a binary quadratic form with integer coefficients and
negative discriminant,
say $g=[\lambda a',\lambda b',\lambda c']$, and $f=[a',b',c']$ with 
$\lambda>0$, then Bernays' result (\ref{R.8}) implies that, as $x$ tends to infinity,
$$B_g(x)\sim C(g){x\over \sqrt{\log x}}.$$
It is easy to see that if $L$ is any arithmetic lattice, then
\begin{equation}
\label{rephraser}
E_L={\sqrt{|D_f|}\over 2}C(f),
\end{equation}
where $f$ is a quadratic form associated to the lattice $L$. Note that if $f$ and $g$ are 
projectively equivalent, then  
$\sqrt{|D_f|}C(f)=\sqrt{|D_g|}C(g)$. We now have:
\begin{Prop}
\label{proptwee}
Let $L$ be a two-dimensional lattice.
The assertion $E_L>E_{\Sigma}$ is equivalent with the assertion that
the minimal value of $\sqrt{|D_f|}C(f)/2$, as $f$ ranges over the primitive
binary quadratic forms of negative discriminant,
is assumed for $f=X^2+XY+Y^2$.
\end{Prop} 

\begin{proof} By Proposition \ref{prop1} we can 
restrict ourselves to arithmetic lattices. The quadratic form associated to an
arithmetic lattice is projectively equivalent with 
a primitive positive definite binary quadratic form of negative discriminant. Vice versa, to a 
quadratic form having 
integer coefficients there corresponds an arithmetic lattice. 
The proof is then completed on invoking (\ref{rephraser}) and noting that
$X^2+XY+Y^2$ is the primitive binary quadratic form associated to $\Sigma$. \qed
\end{proof}

\subsection{On computing the population fraction}
Proposition \ref{proptwee} `reduces' our geometric problem to a problem in number theory, 
namely that of computing $C(f)$. We now discuss some historic results which are related to the explicit evaluation 
of $C(f)$ due to Bernays.

A nonsquare integer $D$ with $D\equiv 0$ or $1({\rm mod~}4)$ is
called a {\it discriminant}. The conductor of the discriminant $D$ is
the largest positive integer $f$ such that $d_0:=D/f^2$ is a
discriminant. If $f=1$, then $D$ is said to be a {\it  fundamental 
discriminant}.
James \cite{James} proved that the number $B_D(x)$ of positive integers $n\le x$ which are 
coprime to $D$ and which are represented by some primitive integral form of discriminant
$D\le -3$ satisfies
$$B_D(x)=J(D){x\over \sqrt{\log x}}+O\Bigl({x\over \log x}\Bigr),$$
where $J(D)$ is the positive constant given by
\begin{equation}
\label{jeedee}
\pi J(D)^2={\varphi(|D|)\over |D|}L(1,\chi_D)\prod_{({D\over p})=-1}{1\over 1-{1\over p^2}},
\end{equation}
and $p$ runs over all primes such that $({D\over p})=-1$. 
Here and in the remainder of the paper implicit constants depend at most on the discriminant
$D$. 
 
 \indent Just as for the characteristic function of $X^2+Y^2$, the characteristic function corresponding
to integers counted for some $x$ by $B_D(x)$ is multiplicative. In both cases the associated
Dirichlet series are very similar and this allowed James to essentially mimic Landau's 
original proof. In 1975 Williams \cite{Williams} reproved James' result in a more
elementary way (essentially along the lines of Rieger \cite{Rieger}, who gave a more
elementary proof of Landau's result). 
However, this reproof only gives a weaker error term.
We like to point out that an even easier proof (but with an even weaker error term) can be
obtained on invoking the following classical result of Wirsing \cite{Wirsing}.

\begin{Thm}
\label{wirsie}
Suppose that $f(n)$ is a multiplicative function such that $f(n)\ge 0$, for $n\ge 1$, and such
that there are constants $\gamma_1$ and $\gamma_2$, with $\gamma_2<2$, such that for every
prime $p$ and for every $\nu\ge 2$, $f(p^{\nu})\le \gamma_1 \gamma_2^{\nu}$. Assume that
as $x\rightarrow \infty$, 
$$\sum_{p\le x}f(p)\sim \tau {x\over \log x},$$
where $\tau>0$ is a constant. Then as $x$ tends to infinity we have
$$\sum_{n\le x}f(n)\sim {e^{-\gamma \tau}\over \Gamma(\tau)}{x\over \log x}
\prod_{p\le x}\left(1+{f(p)\over p}+{f(p^2)\over p^2}+\cdots \right),$$
where $\gamma$ is Euler's constant and $\Gamma(\tau)$ denotes the gamma-function.
\end{Thm}

Let $\xi_D$ be the multiplicative function defined as follows:

\begin{equation*}
\xi_D(p^e)=
\begin{cases}
1 & \text{if $({D\over p})=1;$} \\
1 & \text{if $({D\over p})=-1$ and $2|e;$} \\
0 & \text{otherwise.}
\end{cases}
\end{equation*}

Let $n$ be any integer coprime to $D$. Then $\xi_D(n)=1$ iff $n$ is represented by some
primitive positive integral binary quadratic form of discriminant $D$. It follows that
$B_D(x)=\sum_{n\le x}\xi_D(n)$. It is a consequence of the law of quadratic reciprocity and
the prime number theorem for arithmetic progressions that

\begin{equation}
\label{AP}
\sum_{p\le x}\xi_D(p)=\sum_{p\le x\atop ({D\over p})=1}1\sim {x\over 2\log x}.
\end{equation}

\noindent Thus the conditions of Wirsing's theorem are satisfied and we find that
$$B_D(x)\sim {e^{-\gamma/2}\over \Gamma(1/2)}{x\over \log x}\prod_{p\le x\atop ({D\over p})=1}
{1\over 1-{1\over p}}\prod_{p\le x\atop ({D\over p})=-1}{1\over 1-{1\over p^2}}.$$
By (6) of \cite{Williams} we have the following estimate:


\begin{center}
$\displaystyle \prod_{p\le x\atop ({D\over p}) =1}\Big(1-{1\over p}\Big) = e^{-\gamma/2}\prod_{p|D}\Big(1-{1\over p}\Big)^{-1/2}
\prod_{({D\over p})=-1}\Big(1-{1\over p^2}\Big)^{-1/2}
{L(1,\chi_D)^{-1/2}\over \sqrt{\log x}}
+O\left({1\over \log^{3/2}x}\right)$.
\end{center}

On combining the latter formulae it then follows that 

$$B_D(x)\sim J(D)\frac{x}{\sqrt{\log x}}.$$ 

\indent Indeed on using standard results from the asymptotic theory of arithmetical functions it
is not difficult to improve on James' result. Estimate (\ref{AP}) can be easily sharpened to
$$\sum_{p\le x}\xi_D(p)={1\over 2}\int_2^x{dt\over \log t}+O_m\Bigl({x\over \log^{m}x}\Bigr),$$
for every $m\ge 0$. This in combination with e.g. \cite[Theorem 6]{MoCa} then shows the truth
of the following result:

\begin{Thm}
We have, for every $k\ge 1$,
$$B_D(x)=J(D){x\over \sqrt{\log x}}+\sum_{j=1}^k c_k{x\over \log ^{j+1/2}x}+
O_k\Biggl({x\over \log^{k+3/2-\epsilon}x}\Biggr),$$
where the constants $c_1,c_2,\ldots$ may depend on $D$.
\end{Thm}

\indent James' counting function is artificial in the sense that one would like to
drop the condition that $n$ be coprime to $D$. This was achieved by
Pall \cite{Pall} who proved that the number $C_D(x)$ of positive integers $n\le x$ which are 
which are represented by some primitive integral form of discriminant
$D\le -3$ satisfies
$$C_D(x)=P(D){x\over \sqrt{\log x}}+O\Bigl({x\over \log x}\Bigl),$$
where $P(D)$, {\it Pall's constant}, is computed as follows.
Let $p$ be a prime dividing $D$. 
Let $p^{\prime}$ denote
the primes which satisfy the following condition: if $p>2$ and $p^2 \mid D$ or $p=2$ and 
$D \equiv 0$ or $4 \bmod 16$. Then 
$$P(D)=b_{0}\prod_{p^{\prime}} \Big(1-\frac{1}{{p^{\prime}}^2}\Big)^{-1} \prod_{}
\Big(1+\frac{1}{p^{2k+1}}\Big),$$
where in the second product $D=p^{2k}D^{\prime}$ where $p^2 \nmid D^{\prime}$, $k \geq 1$, and
$\Bigl(\frac{D^{\prime}}{p}\Bigr) \neq -1$, and
$$b_{0}^2=\frac{2h(D)}{w\sqrt{|D|}}\prod_{q} \Big(1-\frac{1}{q^2}\Big)^{-1} \prod_{p^{\prime}}
\Big(1-\frac{1}{p^{\prime}}\Big) \prod_{p\mid D \atop p\neq p^{\prime}} \Big(1-\frac{1}{p}\Big)^{-1},$$
where $q$ runs over all primes such that $\Bigl(\frac{D}{q}\Bigr)=-1$.\\
\indent Let us compute a specific example. 
If $D=-3$, then 
$$P(-3)^2=b_{0}^2=\frac{1}{3} \cdot \frac{1}{\sqrt{3}} \cdot \alpha 
\cdot \frac{3}{2}=\frac{\alpha}{2\sqrt{3}},$$
where $$\displaystyle \alpha=\prod_{q\equiv 2({\rm mod~}3)} \Big(1-{1\over q^2}\Big)^{-1}.$$

\noindent Thus 
$$P(-3)=\frac{1}{\sqrt{2}}\frac{1}{3^{1/4}} \prod_{q\equiv 2({\rm mod~}3)} \Big(1-{1\over q^2}\Big)^{-1/2}.$$
Using Pall's result and the fact that $h(-3)=1$, it then follows that $E_{\Sigma}$ is as given
in (\ref{basaleongelijkheid}). Pall's result allows us to compute $C(f)$ in case the
order associated to $f$ has class number one.\\
\indent Going beyond Pall's work requires genus theory. 
Let $H(D)$ denote the group of strict equivalence
classes of primitive, positive-definite, integral, binary
quadratic forms of discriminant $D$ under Gaussian composition.
Let $G(D)$ denote the genus group of $H(D)$, that is,
$G(D)=H(D)/H(D)^2$. The order $|G(D)|$ of $G(D)$ is a power
of $2$ so that there exists a non-negative integer $t(D)$ such
that $|G(D)|=2^{t(D)}$. The latter quantity is the number
of classes whose order divides $2$, that is, the number of
ambiguous classes in $H(D)$. The value of $t(D)$ is given
as follows (see \cite{cox} or \cite{sw}):

\begin{equation}
\label{td}
t(D)=
\begin{cases}
\omega(D) & \text{if $D\equiv 0({\rm mod~}32);$} \\
\omega(D)-2 & \text{if $D\equiv 4({\rm mod~}16);$} \\
\omega(D)-1 & \text{otherwise,}
\end{cases}
\end{equation}

\noindent where $\omega(D)$ denotes the number of distinct prime factors
in $D$. For example, if $D=-3\equiv 1 \bmod 4$, then $\omega(D)=1$ 
and so there is one genus of forms of discriminant $-3$. 
Note that if $D$ is fundamental, then $t(D)=\omega(D)-1$.
We say that 
$n$ is represented by the genus $G$ of
$G(D)$ if it is represented by at least one class in $G$.
By $g(n,D)$ we denote the number of genera of discriminant $D$
representing $n$. We now turn to the explicit evaluation of $C(f)$ (see page 59 and 115-116 in \cite{B}) which is due to Bernays. 
Namely, we have the following.

\begin{Thm} {\rm (Bernays' Theorem)}. 
\label{bernieboy}
Let $f$ be a positive definite binary quadratic form having discirminant
$D$. Then

\begin{equation}
\label{R.11}
C(f)={J(D)\over 2^{t(D)}}\sum_{n \mid D^{\infty}}\frac{g(n,D)}{n}, 
 \end{equation} 

\noindent where $n\mid D^{\infty}$ means that $n$ divides some arbitrary power of $D$.
\end{Thm}

\noindent It is a classical fact that if $n$ is
represented by a class of discriminant $D$ and $(n,D)=1$, then
$g(n,D)=1$. It is rather more complicated to determine
the value of $g(n,D)$ in case $(n,D)>1$. This was recently
achieved by Kaplan and Williams in \cite{kw} and Sun and Williams in \cite{sw}.
In \cite{kw} they showed that if $g(n,D)>0$, then
$g(n,D)=2^{t(D)-t(D/m^2)}$, where $m$ is the largest integer
such that $m^2|n$ and $m|f$. Note that 
$m^2$ is the largest square dividing $(n,f^2)$. This result together with Theorem
6.1 of \cite{sw} then yields the following result. Here $\nu_p(n)$ denotes the largest power
of the prime $p$ dividing the nonzero integer $n$. 

\begin{Thm}
\label{sunny}
Let $D$ be a discriminant with conductor $f$, $d_0=D/f^2$ and $n$
a natural number. If $(n,f^2)$ is not a square, or there exists
a prime $p$ such that $\nu_p(n)$ is odd and $({d_0\over p})=-1$, then
$g(n,d)=0$. Suppose $(n,f^2)$ is a square and $({d_0\over p})=0,1$ for
every prime $p$ with $\nu_p(n)$ is odd. Then
$g(n,D)=2^{t(D)-t(D/(n,f^2))}$.
\end{Thm}

Using Theorem \ref{sunny} one can evaluate more explicitly the sum

\begin{equation}
\label{sum}
v(D):=\sum_{n|D^{\infty}}{g(n,D) \over n}.
\end{equation}

By Theorem \ref{sunny} 
we have 

\begin{equation}
\label{sunny2}
v(D)=\sum_{m|f}{2^{t(D)-t(D/m^2)}\over m^2}\sum\nolimits'{1\over n_0},
\end{equation}

\noindent where the dash indicates that the sum is over those $n_0$ dividing $D^{\infty}$ such
that $(n_0,f/m)=1$ and there is no prime $p$ such that $2\nmid \nu_p(n_0)$ and
$({d_0\over p})=-1$. 
Note that if $g(n,D)>0$ we can write, by Theorem \ref{sunny}, $(n,f^2)=m^2$, with 
$m|f$ and thus we have $n=n_0m^2$, where $(n_0,f/m)=1$. Furthermore note that
$2\nmid \nu_p(n)$ iff $2\nmid \nu_p(n_0)$. On evaluating the double sum in (\ref{sunny2}) we
obtain

\begin{equation}
\label{formulier}
v(D)={|D|\over \varphi(|D|)}\prod_{p|D \atop ({d_0\over p})=-1}{1\over 1+1/p}
\sum_{m|f}{2^{t(D)-t(D/m^2)}\over m^2}\prod_{p|f/m}\Bigl(1-{1\over p}\Bigr)\prod_{p|f/m \atop ({d_0\over p})=-1}
\Bigl(1+{1\over p}\Bigr).
\end{equation}

Using (\ref{td}) the sum $v(D)$ can be explicitly computed using this formula. Note that it always is
a positive rational number. Also note that if $D$ is a fundamental discriminant, then
$$v(D)={|D|\over \varphi(|D|)}.$$

\begin{example} Take $D=-1984=-2^6\cdot 31$. There are $2^{t(D)}=4$ 
genera of discriminant $-1984$.
We have $G(-1984)=\{I,A,B,AB\}\cong \Bbb Z_2 \times \Bbb Z_2$,
where
\begin{eqnarray}
I&=&\{[1,0,496],[20,\pm4,25]\},\cr
A&=&\{[4,4,125],[5,\pm 4,100]\},\cr
B&=&\{[16,0,31],[7,\pm 2,71]\},\cr
AB&=&\{[16,16,35],[19,\pm 12, 28]\}. \nonumber
\end{eqnarray}
The divisors $n$ of $D^{\infty}$ such that $g(n,D)>0$ are precisely the numbers of the
form $31^a,4\cdot 31^a, 16\cdot 31^a$ and $64\cdot 31^a \cdot 2^b$, where
$a,b\ge 0$ are arbitrary integers. By Theorem \ref{sunny} 
we have $g(n,D)=1,2,4$ and respectively $4$ for these cases. Indeed, if 
$n=31^a$, then the corresponding genera are $I$ and $B$, depending on whether
$a$ is even or odd. If 
$n=4\cdot 31^a$, then the corresponding genera are $I$ and $A$, and $B$ and
$AB$ depending on whether
$a$ is even or odd. In case $n=16\cdot 31^a$ and $n=64\cdot 31^a\cdot 2^b$ the corresponding
genera are $I,A,B$ and $AB$. For example, if $n=4\cdot 31^{2a+1}$, then 
$n$ is represented by $[16,16,35]$ on taking $x=31^a$ and $y=-2\cdot 31^a$ and thus
is represented by $AB$. It follows that
$$v(D)=\sum_{n|D^{\infty}}{g(n,D)\over n}=\left(1+{2\over 4}+{4\over 16}
+{4\over 64}
\sum_{b=0}^{\infty}{1\over 2^b}\right)\sum_{a=0}^{\infty}{1\over 31^a}={31\over 16}.$$
Note that formula (\ref{formulier}) also yields that $v(D)=31/16$. 
\end{example}

\begin{remark} Fomenko \cite{Fomenko} has given an alternative proof of Bernays'
asymptotic result using the theory of multiplicative functions in which the constant $C(f)$ is explicitly computed in case 
$D$ is a fundamental discriminant. Namely, we have (see \cite[Theorem 4]{Fomenko})
for a fundamental discriminant $D$, 
\begin{equation}
\label{fom}
B_{f}(x) \sim {P(D)\over 2^{t(D)}}{x\over \sqrt{\log x}},
\end{equation}
\noindent where $P(D)$ is Pall's constant. It might be interesting to compute $C(f)$
for arbitrary discriminant $D$ using Fomenko's approach.
\end{remark}

\begin{remark}
It might also be of some interest to recover $C(f)$ in general following
Iwaniec's approach to the half-dimensional sieve. Using this sieve (see \cite{iwaniec}), the constant
$C(f_{1})$ was verified for $f_1=X^2 + Y^2$.  
\end{remark}

\section{On explicitly computing the Erd{\H o}s number} 
The explicit formula (\ref{formulier}) for $v(D)$ allows one to
explicitly compute the Erd{\H o}s number $E(D)$.
Note that from (\ref{rephraser}), Theorem \ref{bernieboy}, (\ref{jeedee}), and (\ref{sum}) it follows that

\begin{equation}
\label{edexpliciet}
E(D)={v(D)\over 2^{t(D)+1}}\sqrt{L(1,\chi_D)\varphi(|D|)\over \pi}
\prod_{({D\over p})=-1}\Big(1-{1\over p^2}\Big)^{-1/2},
\end{equation}

\noindent where $v(D)$ is explicitly given by (\ref{formulier}). The latter formula unfortunately does not allow one
to compute $E(D)$ with more that a few decimals of accuracy. A problem
in doing is that the Euler product involved
on direct evaluation (by multiplying consecutive terms
together) can be evaluated with roughly six digit precision
only. However, it turns out that it is possible to express these Euler products in
terms of $L$-series evaluated at integer arguments. To this end note that for $\Re(s)>1/2$, 

\begin{equation}
\label{trick}
\prod_{({D\over q})=-1}
(1-q^{-2s})^{-2} =\frac{\zeta(2s)}{L(2s,\chi_{D})} \prod_{({D\over q})=0} (1-q^{-2s})  \prod_{({D\over q})=-1}
(1-q^{-4s})^{-1}.
\end{equation}

\noindent By recursion we then find from (\ref{edexpliciet}) and (\ref{trick}) the
following formula:

\begin{equation}
\label{trick1}
E(D)={v(D)\over 2^{t(D)+1}}\sqrt{L(1,\chi_D)\varphi(|D|)\over \pi}
\prod_{n=1}^{\infty}\Bigl({\zeta(2^n)\over 
L(2^n,\chi_D)}\prod_{({D\over q})=0} (1-q^{-2^{n}}) \Bigr)^{1/2^{n+1}}.
\end{equation}

\noindent This  
approach was already known to Ramanujan \cite[pp. 60--66]{Berndt} and, 
independently, Shanks \cite[p. 78]{Shanks}.
It can
also be used to deal with more elementary Euler products of the form
$\prod_{p>p_0}(1-f(p)/g(p))$, where $f$ and $g$ are polynomials such that
deg$(f)+2\le {\rm deg}(g)$, see e.g. \cite{Artinconstant}. In the latter
case only values of $\zeta(s)$ at integers are required.

We note that in case $D$ is a fundamental discriminant 
$v(D)=|D|/\varphi(|D|)$ and $t(D)=\omega(D)-1$ and hence

\begin{equation}
\label{trick2}
E(D)={|D|\over 2^{\omega(D)}}\sqrt{L(1,\chi_D)\over \pi\varphi(|D|)}
\prod_{n=1}^{\infty}\Bigl({\zeta(2^n)\over 
L(2^n,\chi_D)}\prod_{({D\over q})=0} (1-q^{-2^{n}}) \Bigr)^{1/2^{n+1}}.
\end{equation}

\section{Some computations of Shanks and Schmid revisited}
We demonstrate our above approach in computing the Erd{\H o}s number (and
hence by (\ref{rephraser}) the Bernays constant $C(f)$), by recomputing
the entries in Table 1 from a
paper by Shanks and Schmid \cite{SS}. They put $C(X^2+nY^2)=b_n$ and
we will follow their notation. The second column in the following table corresponds to the values of $b_{n}$
as computed in \cite{SS} to nine decimal places (for $n=11$, $n=13$ and $n=14$,
approximate values of $b_{n}$ were given). The third column in the table is the computation of $b_{n}$ using 
(\ref{rephraser}) and (\ref{trick1}). 

\newpage

\begin{center}
\begin{tabular}{|c|c|c|}
\multicolumn{3}{c}{} \\ \hline
$n$ & $b_{n}$ & $b_{n}$ \\ \hline
1 & 0.764223654 & 0.7642236535892206629906987311 \\
2 & 0.872887558 & 0.8728875581309146129200636834 \\
3 & 0.638909405 &  0.6389094054453438822549426747 \\
4 & 0.573167740 & 0.5731677401919154972430240483 \\
5 & 0.535179999 & 0.5351799988649545413027199090 \\
6 & 0.558357114 & 0.5583571140895246274460701041 \\
7 & 0.543539641 & 0.5435396411014846926771211300 \\
8 & 0.436443779 & 0.4364437790654573064600318417 \\
9 & 0.424568696 & 0.4245686964384559238837215172 \\
10 & 0.473558100 & 0.4735580999381557098419651553 \\
11 & $\approx$ 0.677 & 0.6773880181341740551427831009 \\
12 & 0.399318378 & 0.3993183784033399264093391717 \\
13 & $\approx$ 0.420 & 0.4207205175783009914997595500 \\
14 & $\approx$ 0.563 & 0.5634867715862649042931719141 \\
16 & 0.334347848 & 0.3343478484452840400584306948 \\
20 & 0.401384999 & 0.4013849991487159059770399317 \\
24 & 0.279178557 & 0.2791785570447623137230350520 \\
27 & 0.496929538 & 0.4969295375686007973093998581 \\
64 & 0.274642876 & 0.2746428755086261757622823564 \\
96 & 0.209383918 & 0.2093839177835717352922762890 \\
256 & 0.259716632 & 0.2597166322744617096882452719 \\ \hline
\end{tabular}
\end{center}

\section{Proof of Theorem 1}
 
The idea of the proof is to use a lower bound estimate for $\varphi(|D|)$ combined with an 
upper bound estimate for $\omega(D)$ to show that $E(D) > E(-3)$ for all $|D| \geq D_{0}$,
with
$D_{0}$ an explicit number. In the range $|D|<D_{0}$ one then determines
those $D$ for which the quickly computed lower bound for $E(D)^2$ given
in (\ref{globbie2}) does not exceed $E(-3)^2$.  
For these values of $D$ one then computes $E(D)$ using (\ref{trick1})
and compares with $E(-3)$. We now prove Theorem 1.

\begin{proof}  Note 
that $h(D)\geq 2^{t(D)}$, $v(D) \geq 1$ and that the Euler product in (\ref{edexpliciet})
exceeds one. Using these 
trivial lower bounds and (\ref{edexpliciet}) we infer that

\begin{equation}
\label{globbie}
E(D) \geq \Bigg( \frac{1}{2^{t(D)+1}w(D)} \frac{\varphi(|D|)}{|D|^{1/2}} \Bigg)^{1/2},
\end{equation}

\noindent where we used that $L(1,\chi_D)=2\pi h(D)/(w(D)\sqrt{|D|})$. It is well-known that
in case $D$ is a fundamental discriminant $w(D)\ne 2$ if and only if $D=-3$ or $D=-4$.
Using the observation that that order for the discriminant $D$ is the $\mathbb Z$-module generated
by 1 and $f(D+\sqrt{D})/2$
(cf. \cite[Lemma 7.2]{cox}), where $f$ is the conductor, one sees that $w(D)=2$ unless $D=-4$ or $D=-3$.
In the rest of the proof we assume that $|D|\ge 5$. Then 

\begin{equation}
\label{globbie2}
E(D)^2 \geq \frac{\varphi(|D|)}{2^{t(D)+2}\sqrt{|D|}}. 
\end{equation}

Put $g(n)=\varphi(n)/(2^{\omega(n)}\sqrt{n})$. Note that $g$ is a multplicative
function of $n$. If $n=\prod_{i=1}^m q_i^{e_i}$ denotes the canonical factorisation of $n$, then

$$g(n)=\prod_{i=1}^m {1\over 2}q_i^{e_i/2-1}(q_i-1)\ge \prod_{i=1}^m{1\over 2}\Bigl(\sqrt{q_i}-{1\over \sqrt{q_i}}\Bigr).$$
We let $p_1=2,~p_2=3,\cdots$ denote the consecutive
primes. Note that $\sqrt{x}-{1\over \sqrt{x}}$ is strictly increasing with $x$.
It thus follows that
$$g(n)\ge \prod_{i=1}^m {1\over 2}\Bigl(\sqrt{p_i}-{1\over \sqrt{p_i}}\Bigr).$$
If $n$ is odd, then we similarly have

\begin{equation}
\label{oneven} 
g(n)\ge \prod_{i=2}^{m+1}{1\over 2}\Bigl(\sqrt{p_i}-{1\over \sqrt{p_i}}\Bigr).
\end{equation}

\noindent From (\ref{td}) and (\ref{globbie2}) one infers that

\begin{equation}
\label{aadd} 
E(D)^2\ge \alpha(D) g(D_{\rm odd}),
\end{equation} 

\noindent where 
\begin{equation*}
\alpha(D)=
\begin{cases}
1/4 & \text{if $D\equiv 12({\rm mod~}16);$} \\
{1/2\sqrt{2}} & \text{if $D\equiv 8({\rm mod~}16)$ or $D\equiv 0({\rm mod~}32);$} \\
1/2 & \text{if $D\equiv 1({\rm mod~}4)$,~$D\equiv 0({\rm mod~}16)$, or 
$D\equiv 4({\rm mod~}16)$,}
\end{cases}
\end{equation*}

\noindent and $D_{\rm odd}$ denotes the largest odd divisor of $D$.

First assume that $D\equiv 1({\rm mod~}4)$ (thus
$\alpha(D)=1/2$ and $t(D)=\omega(D)-1$). Then, from (\ref{aadd}) 
and (\ref{oneven}) we infer that
$$2E(D)^2\ge \prod_{i=2}^{\omega(D)+1}{1\over 2}\Bigl(\sqrt{p_i}-{1\over \sqrt{p_i}}\Bigr).$$ 
If $\omega(D)>3$ it follows from the latter inequality that $E(D)>0.66>E(-3)$.
So let us assume that $\omega(D)\le 3$. It now follows, using that 

\begin{equation}
\label{nicolas}
{\varphi(n)}>e^{-\gamma}{n\over \log \log n},
\end{equation}

\noindent for
all odd integers $n\ge 17$ (see \cite{CLMS}), that for $|D|\ge 19$ we have
$$E(D)^2\ge {\varphi(|D|)\over 16\sqrt{|D|}}\ge {e^{-\gamma}\over 16}{\sqrt{|D|}\over \log \log |D|}.$$
From this estimate one infers that $E(D)>E(-3)$ for $|D|\ge 217$. For the
$D$ with $D\equiv 1({\rm mod~}4)$ and $7\le |D|\le 215$ one checks that
$$\left({\varphi(|D|)\over 2^{\omega(D)+1}\sqrt{|D|}}\right)^{1/2}>0.6>E(-3),$$
except for $D=-15$. A direct computation shows that $E(-15)=0.9719612\cdots >E(-3)$. 

The remaining cases are dealt with similarly: on noting that the right hand side
of (\ref{oneven}) is monotonically increasing for $m\ge 2$ one uses (\ref{aadd}) to
obtain an upper bound for $\omega(D)$. From this upper bound, (\ref{globbie2})
and (\ref{nicolas}), one then finds an integer $D_0$
such  that if $E(D)>E(-3)$, then $|D|<D_0$. For the
discriminants $D$ with $|D|<D_0$ one then computes the
discriminants $D$ for which the left hand side of (\ref{globbie2}) does
not exceed $E(-3)^2$.
For these $D$ values one then computes $E(D)$ using (\ref{trick1}). 
One finds that for all these values of $D$ one
has $E(D)>E(-3)$. In this way it is seen
that $E(D)$ is minimal for $D=-3$.

To prove the second assertion note that in the above argument one can replace
$E(-3)$ with any real number $r$. In the end one is left with a finite
list of $D$ for which $E(D)<r$. \qed
\end{proof}

\begin{example} If $r=1$, then one finds the following list.

\begin{center}
\begin{tabular}{|c|c|}
\multicolumn{2}{c}{} \\ \hline
$D$ & $E(D)$  \\ \hline
$-3$ & 0.5533117758324795595155817776 \\
$-4$ &  0.7642236535892206629906987311 \\
$ -7$ & 0.9587138120398867707178043483 \\
$-15$ & 0.9719612596359906049817562980 \\ \hline
\end{tabular}
 
\end{center} 
 
\vspace{.1in} 

\noindent Thus the second smallest lattice is given by the maximal order with $D=-4$ (the
square lattice)
and the third and fourth smallest lattices by $D=-7$ and $D=-15$ respectively. 
\end{example}

 \begin{remark} The inequality (\ref{nicolas}) is quite subtle.
Let $N_k=2\cdot 3\cdots p_k$ be the product of the first $k$ primes, then
if the Riemann Hypothesis is true (\ref{nicolas}) is false for every
integer $n$ with $n=N_k$. On the other hand, if the Riemann Hypothesis is
false then there are infinitely many integers $k$ for which $n=N_k$ does
satisfy (\ref{nicolas}). See Nicolas \cite{nic} for a proof of this
interesting result.
\end{remark}

\section*{Acknowledgement} This paper owes much to an inspiring discussion
with Prof.~Don Zagier in which he convinced the first author that proving Theorem \ref{main}
should be doable. The authors would like to thank Valentin Blomer for his
helpful comments regarding Bernays' thesis and K.S.~Williams for making his preprint \cite{sw} available.
It is also a pleasure to thank UCD graduate student Raja Mukherji for his suggestions which greatly improved the efficiency of the GP/PARI program which was used in Sections 4 and 5. Finally, the authors thank the Max-Planck-Institut f{\"u}r 
Mathematik in Bonn for its hospitality and support during the preparation of this paper.

\EMdate{20, Octobre 2006}

\begin{address}
Pieter Moree \\
Max-Planck-Institut f\"ur Mathematik
Vivatsgasse 7
D-53111 
Bonn
Germany
\email{moree@mpim-bonn.mpg.de}
\end{address}

\begin{address}
Robert Osburn \\
School of Mathematical Sciences
University College Dublin
Belfield
Dublin 4 
Ireland
\email{robert.osburn@ucd.ie}
\end{address}  

\end{document}